\documentclass[conference,letterpaper]{IEEEtran}
\usepackage{graphicx}
\usepackage{subfigure,stfloats}
\usepackage{amssymb}
\usepackage{tikz}

\usepackage[utf8]{inputenc} 
\usepackage[T1]{fontenc}
\usepackage{url}
\usepackage{ifthen}
\usepackage{cite}
\usepackage[cmex10]{amsmath} % Use the [cmex10] option to ensure complicance
\newtheorem{theorem}{Theorem}
\newtheorem{lemma}{Lemma}
\newtheorem{corollary}{Corollary}
\newtheorem{proposition}{Proposition}
\newtheorem{problem}{Problem}
\newtheorem{remark}{Remark}

\newcommand\SSSS[1][n]{S^{(4)}_{#1}}
\newcommand\SSS[1][n]{S^{(3)}_{#1}}
\renewcommand\SS[1][n]{S^{(2)}_{#1}}

\begin{document}
\title{On $(2n/3-1)$-Resilient $(n,2)$-Functions} 
\author{%
   \IEEEauthorblockN{Denis~S.~Krotov%
}
   \IEEEauthorblockA{Sobolev Institute of Mathematics,
                     Novosibirsk 630090, Russia\\
                     Email: krotov@math.nsc.ru, dk@ieee.org}
}

\maketitle

\begin{abstract}\boldmath
A $\{00,01,10,11\}$-valued function on the vertices of the $n$-cube
is called a $t$-resilient $(n,2)$-function 
if it has the same number of $00$s, $01$s, $10$s and $11$s
among the vertices of every subcube of dimension $t$.
The Friedman and Fon-Der-Flaass bounds on the correlation immunity order say
that such a function must satisfy $t\le 2n/3-1$; moreover,
the $(2n/3-1)$-resilient $(n,2)$-functions correspond
to the equitable partitions of the $n$-cube with the quotient matrix
$[[0,r,r,r],[r,0,r,r],[r,r,0,r],[r,r,r,0]]$, $r=n/3$.
We suggest constructions of such functions and corresponding partitions,
show connections with Latin hypercubes and binary $1$-perfect codes,
characterize the non-full-rank and the reducible functions from the considered class,
and discuss  the possibility to make a complete characterization of the class.
\end{abstract}

\def\Z{\mathbb Z}
\def\F{\mathbb F}
\def\E{\mathbb E}
\def\C{\mathbb C}
\def\R{\mathbb R}
\def\V{V_{m,n',n''}}
\def\inn#1#2{{\langle #1,#2 \rangle}}
\def\innb#1#2{\inn{\vec{#1}}{\vec{#2}}}
\newcommand\zinn[3][{}]{{[ #2,#3 ]_{#1}}}
\newcommand\zinnb[3][{}]{\zinn[#1]{\vec{#2}}{\vec{#3}}}
\def\vec#1{{\boldsymbol{#1}}}
\def\2{2\cdot }
\newcommand\Tr[1][{}]{\mathrm{Tr}_{#1}}
\newcommand\wt{\mathrm{wt}}
\newcommand\wa{\widetilde{\phantom{wa}}\!\!\!\!\!\!\!\mathrm{wt}}
\def\cC{\overline C}
\def\CcC{(C,\overline C)}

\section{Introduction}
\renewcommand\thefootnote{}
\footnotetext{This work was funded by the Russian Science Foundation (RSF)\linebreak 
under grant  18-11-00136}%

The set $Q_n:=\{0,1\}^n$ of all $n$-words over $\{0,1\}$ is called 
the $n$-cube. This set forms a linear $n$-dimensional space over the field GF$(2)$
with the standard basis $\vec e_1=10...0$, $\vec e_2=010...0$, \ldots, $\vec e_n=0...01$.
The same term $Q_n$ is used to denote the graph on $\{0,1\}^n$,
where two words are adjacent if and only if they differ in exactly
one position. The number of ones in a word $\vec v$ from $Q_n$ is referred to as
the \emph{weight} of $\vec v$, $\wt(\vec v)$.
A \emph{$k$-subcube} of $Q_n$ is a subgraph isomorphic
to $Q_k$ (and also the corresponding subset of vertices).
%!!!!! discuss the definition of a subcube
A (vectorial Boolean) function from $Q_n$ to $Q_m$ 
is called an $(n,m)$-function.
The $(n,1)$-functions correspond to the usual,
non-vectorial, Boolean functions.
An $(n,m)$-function is called \emph{balanced} 
if it possesses each of the $2^m$ values exactly 
$2^{n-m}$ times.
An function from $Q_n$ to some finite set $S$ is called 
\emph{correlation immune}
of order $t$ if the proportion of occurrences of the values is the same
in all $l$-subcubes of $Q_n$ with $l\ge n-t$. So, 
$$ |f^{-1}(s) \cap Q'_l |= |f^{-1}(s)| / 2^{n-l}$$
for every $s\in S$ and every $l$-subcube $Q'_l$, $l\ge n-t$.
An $(n,m)$-function is called \emph{$t$-resilient} if it is both balanced 
and correlation immune of order $t$. Equivalently, if
for every $\vec y \in Q_m$, the characteristic function
\begin{IEEEeqnarray}{c}\label{eq:fy}
 f_{\vec y}(\vec x):=
\begin{cases}
1 &\mbox{if } f(\vec x)=\vec y\cr
0 &\mbox{otherwise}
\end{cases}
\end{IEEEeqnarray}
of $f^{-1}(\vec y)$ is correlation immune of order $t$ and has exactly $2^{n-m}$ ones.
Resilient functions play important role in cryptography,
see e.g. \cite{Carlet:vect}.

In \cite{Friedman:92}, Friedman derived the bound 
\begin{IEEEeqnarray}{c}\label{eq:fb}
 \frac{n-t-1}{n} \ge \frac{2^{m-1}-1}{2^{m}-1}
\end{IEEEeqnarray}
 for $t$-resilient $(n,m)$-functions or, more generally,
for $(n,1)$-functions of correlation-immunity order $t$ with $2^{n-m}$ ones. 
For $m=2$, this bound turns to 
\begin{IEEEeqnarray}{c}\label{eq:b}
 t \le \frac{2n}{3}-1, 
\end{IEEEeqnarray}
and in the current thesis we are interested 
in the $t$-resilient $(n,2)$ functions that meet \eqref{eq:b}
with equality.
Fon-Der-Flaass \cite{FDF:CorrImmBound} proved that 
the correlation-immunity order $t$ of
any non-constant unbalanced Boolean $(n,1)$ function
meets \eqref{eq:b}.
(Note that there are interesting classes of unbalanced Boolean $(n,1)$ functions
with more than $2^{n-2}$ and less than $2^{n-1}$ 
ones that attend the Fon-Der-Flaass
bound, see \cite{Tarannikov2000}, \cite{FDF:12cube.en}.) 
Since for $m\ge 2$ the functions $f_{\vec y}$, see \eqref{eq:fy},
are not balanced, we see that the class of functions we study 
lies on  both the Friedman bound and the Fon-Der-Flaass bound.
It occurs that attending any of these two bounds implies that the function
belongs to the class of very regular objects, known as equitable partitions.

A partition $(C_i)_{i\in I}$ (where $I$ is a finite index set) 
of the vertices of a graph $\Gamma$ 
is called an \emph{equitable partition}
(perfect coloring, regular partition, partition design)
with a quotient matrix $[[S_{i,j}]]_{i,j\in I}$
if for every $i,j\in I$ 
and for every $c\in C_i$ 
one has $|\Gamma(c)\cap C_j|=S_{i,j}$,
where $\Gamma(c)$ is the neighborhood of $c$ in $\Gamma$.

As was proved by Fon-Der-Flaass \cite{FDF:CorrImmBound}, any non-constant unbalanced 
Boolean function attaining the bound \eqref{eq:b} on the order $t$ of correlation immunity
corresponds to an equitable partition $(f^{-1}(1),f^{-1}(0))$ with a quotient matrix
$[[a,b],[c,d]]$, where $a+b=c+d=n$ and $a-c=d-b=-n/3$. 
On the other hand, Potapov \cite{Potapov:2010,Pot:2012:color} proved
that any Boolean  function with $2^{n-m}$ ones whose correlation-immunity order $t$ satisfies \eqref{eq:fb} with equality corresponds to an equitable partition $(f^{-1}(1),f^{-1}(0))$ with the quotient matrix
$[[0,n],[n/(2^m-1),n-n/(2^m-1)]]$.
From any of these results, we get the following.
\begin{proposition}\label{p:eqp}
For any $(2n/3-1)$-resilient $(n,2)$-function $f$, 
the partition $(f^{-1}(00),f^{-1}(01),f^{-1}(10),f^{-1}(11))$ of $Q_n$
is equitable with the quotient matrix
\begin{IEEEeqnarray}{c}\label{eq:S4}
\SSSS:=
\left[
\begin{array}{cccc}
 0 & r & r & r \\
 r & 0 & r & r \\
 r & r & 0 & r \\
 r & r & r & 0 
\end{array}
\right], \qquad r=\frac n3
\end{IEEEeqnarray}
\begin{remark}
 In general, one can see from \cite{Potapov:2010} that 
any $t$-resilient $(n,m)$-function such that 
\eqref{eq:fb} holds with equality induces an equitable partition
of $Q_n$ into $2^m$ cells with the following quotient matrix:
the diagonal elements are $0$, while every non-diagonal element is $n/(2^m-1)$.
\end{remark}

In this thesis, we study the $(2n/3-1)$-resilient $(n,2)$-functions or, eqiuvalently,
the equitable $4$-partitions of $Q_n$ with the quotient matrix \eqref{eq:S4}, and 
the related (as we see from \eqref{p:eqp}) class of  order-$(2n/3-1)$ correlation-immune $(n,1)$-functions
with $2^{n-2}$ ones, or, eqiuvalently, the equitable $2$-partitions of $Q_n$ with  the quotient matrix 
\begin{IEEEeqnarray}{c}\label{eq:S2}
\SS:=
\left[
\begin{array}{cc}
        0 &   3r   \\
 r & 2r
\end{array}
\right],\qquad r=\frac n3.
\end{IEEEeqnarray}
\end{proposition}
We will refer to the equitable partitions of $Q_n$ with the quotient matrix $\SSSS$ or $\SS$ as the 
$\SSSS$-partitions or $\SS$-partitions, respectively.

In the next section,
we discuss relations between the classes of 
$\SS$-partitions and $\SSSS$-partitions and
announce rhe results of computational classification of 
$5$-resilient $(9,2)$-functions (Theorem~\ref{th:N9}).
Section~\ref{s:semi} is devoted to the concept of rank of $(n,2)$-functions
and related partitions. The $\SS$-partitions are characterized it terms
of multifold $1$-perfect codes (Theorem~\ref{th:semi}), giving a powerful construction 
of such partitions and $(2n/3-1)$-resilient $(n,2)$-functions.
In Section~\ref{s:conc}, we consider connections of the class of 
$(2n/3-1)$-resilient $(n,2)$-functions with the class of 
Latin $n/3$-cubes of order $4$; 
a related concatenation construction presented in Theorem~\ref{th:conc}.
In Section~\ref{s:contr}, we discuss the existence of 
$(2n/3-1)$-resilient $(n,2)$-functions that cannon be constructed
by any way suggested in the previous sections.
Section~\ref{s:conc} contains the conclusion.

%=============================================
%=============================================
%=============================================
\section{$\SSSS$-partitions vs $\SS$-partitions}\label{s:42}

Given an $\SSSS$-partition $(C_{00},C_{01},C_{10},C_{11})$, each of its cells,
together with the complement, forms a  $\SS$-partition. 

\begin{problem}\label{probl:2to4}
Is the first cell of an $\SS$-partition always occurs as a cell of some 
$\SSSS$-partition?
\end{problem}

In the next sections, we will see that the answer is positive for some special classes
of $\SS$-partition (semilinear $\SS$-partitions and reducible $\SS$-partitions). 

\begin{lemma}\label{l:23}
The second cell $\cC$ of any $\SS$-partition $(C,\cC)$ can be split, $\cC=C_1\cup C_2$,
to form an equitable partition $(C,C_1,C_2)$ of $Q_n$ with the quotient matrix 
\begin{IEEEeqnarray}{c}\label{eq:S3}
\SSS:=
\left[
\begin{array}{ccc}
 0 & r & 2r \\
 r & 0 & 2r \\
 r & r & r 
\end{array}
\right],\qquad r=\frac n3.
\end{IEEEeqnarray}
\end{lemma}
Here and below, we use the notation $\overline C$ for the complement of $C$.
\begin{IEEEproof}
 Take $C_1:=C+0...01$. The set $C$ is independent,
 hence $C\cap C_1=\emptyset$. Therefore,
 the equitability (and the quotient matrix) of $(C,C_1,\cC \backslash C_1)$ is 
 straightforward from that of $(C,\cC)$ and $(C_1, C\cup \cC\backslash C_1)$.
\end{IEEEproof}

\begin{problem}\label{probl:3to4}
Are the first and second cells of an $\SSS$-partition always occur as cells of some 
$\SSSS$-partition?
\end{problem}

A positive answer to   Problem~\ref{probl:3to4} would imply a positive
answer  to   Problem~\ref{probl:2to4}.

Another connection of $4$- and $2$-partitions from the considered classes 
is the following fact.

\begin{lemma}\label{l:42}
 The $\SSSS$-partitions are in one-to-one correspondence 
 with the $\SS[n+3]$-partitions $(C,\overline C)$ such that $C=C+0...0111$.
\end{lemma}
\begin{IEEEproof}[A sketch of proof]
 Given an $\SSSS$-partition $(C_0,C_1,C_2,C_3)$, define 
 $c_0=000$, $c_1=001$, $c_2=010$, $c_3=010$,
 $c_4=100$, $c_5=101$, $c_6=110$, $c_7=110$,
 \begin{IEEEeqnarray}{c}\label{eq:4to2}
  C = \sum_{i=0,1,2,3} C_i \times \{c_i,c_{i+4}\}.
 \end{IEEEeqnarray}
 The required properties of the partition $(C,\overline C)$
 are straightforward to check.
Conversely, given an $\SS$-partition $(C,\overline C)$, define 
$C_0$, \ldots, $C_7$ such that 
 $$
  C = \sum_{i=0}^7 C_i \times \{c_i\}.
 $$
 If $C=C+0...0111$, then $C_i=C_{i+4}$, $i=0,1,2,3$, and we have \eqref{eq:4to2}. 
 The required properties of $(C_0,C_1,C_2,C_3)$
 are straightforward.
\end{IEEEproof}

In \cite{KroVor:CorIm}, 
among other results, it was established that
there are exactly $16$ isomorphism classes of $\SS[12]$-partitions.
Eight of them, after some coordinate permutation, 
are invariant under translation by $0...0111$.
However, different coordinate permutations 
can result in nonisomorphic $\SSSS[9]$-colorings,
in the mean of Lemma~\ref{l:42}. 
In the ordering of \cite{KroVor:CorIm}, 
nonisomorphic $\SS[12]$-partitions
correspond to $1$, $1$, $2$, $1$, $1$, $2$, $0$, $1$, $0$,  
$0$,  $0$,  $0$,  $0$,  $0$, $1$, $0$
nonisomorphic $\SSSS[9]$-partitions, respectively. The total number of the isomorphism classes is $10$:
\begin{theorem}\label{th:N9}
 There are exactly $10$ non-equivalent $5$-resilient $(9,2)$-functions,
 where two functions $f$ and $g$ are \emph{equivalent}
 if $f(\vec x)\equiv \beta(g(\alpha(\vec x)))$
 for some isomorphism (isometry) $\alpha$ of $Q_n$
 and some permutation $\beta$ of $Q_2$.
 Among these $10$ functions, exactly one is linear and exactly one
 is full rank; the remaining $8$ are semilinear (see the definitions in the next section).
\end{theorem}

%=============================================
%=============================================
%=============================================
\section{Non-full-rank partitions, connection with multifold perfect codes}\label{s:semi}
In this section we characterize the $\SS$-partitions of deficient rank in terms
of equitable partitions of $Q_m$, 
$m=2n/3$ with other parameters,
which are closely related with the theory of $1$-perfect codes.

%=============================================
%=============================================
\subsection{The non-full-rank $\SS$-partitions}\label{s:nfr}

The (affine) \emph{rank} of a set $C\subset Q_n$
is the dimension of its affine span.
By the rank of an $\SS$-partition $(C,\overline C)$, 
we will mean the rank of $C$.
As $|C|=2^{n-2}$, the rank can be only
$n-2$ if $C$ is linear,
or $n-1$, in which case we say that $C$ is \emph{strictly semilinear},
or $N$, if $C$ is \emph{full rank}. 
We say that $C$ is \emph{semilinear},
or \emph{non-full-rank}, if it is either linear of strictly semilinear.
The \emph{dual} $C^\perp$ of $C$ is defined as the set of all words
$\vec v$ in $Q_n$ such that $\langle \vec c, \vec v\rangle=0$ 
for all $\vec c$ in $C$ or $\langle \vec c, \vec v\rangle=1$ 
for all $\vec c$ in $C$, where 
$$ \langle (c_1,...,c_n),(v_1,...,v_n)\rangle := c_1v_1+...+c_n v_n. $$
The following fact is well known.
\begin{lemma}\label{l:du}
If $\CcC$ is an equitable partition of $Q_n$ with the quotient matrix $[a,b],[c,d]]$, 
then the nonzero words of $C^\perp$ have weight $(b+c)/2$.
In particular, for the $\SS$-partitions, this weight is $2n/3$. 
\end{lemma}
It immediately follows that there is only one linear $\SS$-partition,
up to isomorphism. For a semilinear  $\SS$-partition,
there is only exactly one nonzero dual word. W.l.o.g,
we can assume that all its zeros are in the last $n/3$ positions.

\begin{remark}
 We defined the rank and the dual of a code in the affine sence, 
 to be invariant with respect to the translations of the space.
 In general, this approach is convenient while considering nonlinear
 vector sets whose main properties are invariant with respect to
 the space isometries. The readers who prefer 
 using the classical (linear-sence)
 concepts of the rank and the dual can treat the material below assuming
 that $C$ always contains the all-zero word.
 \end{remark}

\begin{theorem}\label{th:semi}
 The $\SS$-partitions $\CcC$ such that 
 $\vec v:=\underbrace{1\ldots1}_{2n/3}\underbrace{0\ldots0}_{n/3}\in C^\perp$
 are in one-to-one correspondence with the equitable
 partitions of $Q_{2n/3}$ with the quotient matrix 
 \begin{IEEEeqnarray}{c}\label{ex:eperf}
 \left[
 \begin{array}{ccc}
  0 & 0 & 2r \\
  0 & 0 & 2r \\
  r & r & 0
 \end{array}
 \right],
 \qquad
 r=\frac{n}3.
 \end{IEEEeqnarray}
\end{theorem}
\begin{IEEEproof}[A sketch of proof]
 Consider the case $\langle \vec c, \vec v\rangle=0$ for all $\vec c \in C$.
 Define $C' = \{ c \mid \langle \vec c, \vec v\rangle=0 \} \backslash C$
 and
 $C'' = \{ c \mid \langle \vec c, \vec v\rangle=1 \}$.
 Since $(C\cup C', C'')$ is an equitable partition with the quotient matrix
 $$ \left[
 \begin{array}{cc}
  r & 2r \\
  2r & r  
 \end{array}
 \right],
 $$
 it is easy to check that $(C,C',C'')$ is an $\SSS$-partition.
 
 Next, we define 
 $D=\{ \vec x \mid \vec x\underbrace{0\ldots0}_{n/3}\in C\}$,
 and similar $D'$ and $D''$. 
 It happens that  $D \cup D'=Q_{r,\mathrm{ev}}$ 
 and $D''=\overline Q_{r,\mathrm{ev}}$,
 where $Q_{r,\mathrm{ev}}$ is the set of even-weight words of $Q_r$.

 (i) $C''=\{\vec x \vec y \mid \vec x \in D'', \ \vec y \in Q_r\}$, which is straightforward from the definitions
 of $C''$ and $D''$;
 
 (ii)  $C=\{\vec x \vec y \mid \vec x \in D, \ \vec y \in Q_{r,\mathrm{ev}}\} \cup \{\vec x \vec y \mid \vec x \in D', \ \vec y \in Q_{r,\mathrm{od}}\} $,
 which comes from the definition of $D$, the cardinality of $C$, $|C=2^n/4$, and the fact that $C$ is an independent set;
 
 (iii)  $C'=\{\vec x \vec y \mid \vec x \in D', \ \vec y \in Q_{r,\mathrm{ev}}\} \cup  \{\vec x \vec y \mid \vec x \in D, \ \vec y \in Q_{r,\mathrm{od}}\} $,
 similarly.

 Now, since a vertex $(\vec x,10...0)$, $\vec x \in D''$, has $r$ neighbors from $C$,
 we see that $\vec x$ has $r$ neighbors from $D$. Immediately, we get that $(D,D',D'')$
 is an equitable partition with the quotient matrix \eqref{ex:eperf}.
 
 Conversely, having an equitable partition $(D,D',D'')$ 
 with the quotient matrix \eqref{ex:eperf} and defining $C$, $C'$, $C''$ by (i), (ii), (iii),
 we get an $\SS$-partition. 
 If $C''=\overline Q_{r,\mathrm{ev}}$, 
 then $\langle \vec c, \vec v\rangle=0$ for all $\vec c \in C$.
 If $C''=Q_{r,\mathrm{ev}}$, 
 then $\langle \vec c, \vec v\rangle=1$ for all $\vec c \in C$.
\end{IEEEproof}

From the considerations abobe, it is not difficult to conclude the following.

\begin{corollary}\label{c:semi1}
 The answer to Problem~\ref{probl:2to4} is positive for the class of semilinear
 $\SS$-partitions.
\end{corollary}

%=============================================
%=============================================
\subsection{Multifold perfect codes}\label{s:perf}

In this subsection, we discuss equitable partitions with the quotient matrix \eqref{ex:eperf}
and their connections with perfect codes. A \emph{$1$-perfect code} in $Q_m$ is a set $C$ such that 
$|C\cap B|=1$ for every ball $B$ or radius $1$. More general, a \emph{$t$-fold $1$-perfect code} in $Q_m$ is a set $C$ such that 
$|C\cap B|=1$ for every ball $B$ or radius $1$.
Equivalently, $(C,\overline{C})$
is an equitable partition with quotient
matrix $[[t-1,m-t+1],[t,m-t]]$.
After adding the all-parity-check bit to all words of $C$, we obtain
a \emph{$t$-fold extended $1$-perfect code} $C'$ such that
$(C', Q_{m+1,\mathrm{ev}} \backslash C' ,\overline Q_{m+1,\mathrm{ev}})$
is an equitable partition with the quotient
matrix  \eqref{ex:eperf}, $r=(m+1)/2$.
So, Theorem~\ref{th:N9} establishes a connection between 
the non-full-rank $\SS$-partitions and
the $r$-fold $1$-perfect codes of length $2r-1$, $r=n/3$.
$t$-fold $1$-perfect codes can be easily constructed 
constructed from an arbitrary usual ($1$-fold) $1$-perfect code $C$,
as the union of $t$ translations of $C$ by vectors of weight $1$
(it should be noted that not all $t$-fold codes can be treated in such a way \cite{KroPot:nonsplittable}).
The binary $1$-perfect codes have been extensively studied since the first construction of nonlinear such
codes by Vasil'ev \cite{Vas:nongroup_perfect.en}, but currently,
even the asymptotic of the log\,log of the number of such codes is not known.
It is hard to expect a constructive characterisation of the class
of binary $1$-perfect codes.
The problem of the characterization of all
$r$-fold $1$-perfect codes of length $2r-1$
is not formally harder or simpler of the same problem
for $1$-fold $1$-perfect codes:
we cannot construct all $r$-foll $1$-perfect codes from all $1$-perfect or vice versa.
But intuitively, the problems belong
to the same class: in both cases, we have classes
of equitable $2$ partitions with doubly-exponential (in length)
grows of the number of objects.
The main problem in this area is the 
asymptotic of the log\,log of the number of objects 
(note that this asymptotic is the same
for the number of different object
and the number of nonisomorphic objects),
which is known bounded 
by $n/2$ from below (from \cite{Vas:nongroup_perfect.en}
and similar approach for $t>1$) and by $n$ from above (a trivial bound),
for $t$-perfect binary codes of length $n$,
for any $t$.

%=============================================
%=============================================
\subsection{The rank of $\SSSS$-partitions}\label{s:rank4}

It remains to discuss the rank of an $\SSSS$-partition,
or, equivalently,
a $(2n/3-1)$-resilient $(n,2)$-function. Given such function $f$, 
its (affine) rank
can be defined as the dimension of the affine span if its graph
$$G(f):=\{\vec x\vec y \mid \vec x \in Q_n, \vec y=f(\vec x)\in Q_2\}$$ 
and is obviously 
one of $n$, $n+1$, $n+2$.
A nice property of $Q_2$ is that any permutation $\pi$ of its element is an affine transform;
so, the $(n,2)$ function $\pi(f(\cdot))$ has the same rank as $f(\cdot)$.
This means that we can well define the rank of an $\SSSS$-partition
as the rank of the associated $(n,2)$-function, without caring about the order of the partition cells. Readily, the rank of an $\SSSS$-partition is connected with the rank of the associated 
$\SS$-partition, in the sense of Lemma~\ref{l:42}.
\begin{lemma}\label{l:rank42}
 Let $\vec C$ be an $\SSSS$-partition,
 and let $(C,\cC)$ be the associated $\SS$-partition, in the sense of Lemma~\ref{l:42}.
 Then $\mathrm{rank}(\vec C)+1=\mathrm{rank}(C)$.
\end{lemma}
\begin{IEEEproof}
 Straightforward from
 $C=\{\vec x0 \mid \vec x \in G(f)\}+0...0111$, where $f$ is the $(n,2)$-function corresponding
 to $\vec C$.
\end{IEEEproof}
In this way, we see that the non-full-rank $\SSSS$-partitions are also related to
multifold $1$-perfect codes.

%=============================================
%=============================================
%=============================================
\section{Connection with latin hypercubes of order $4$, 
reducible $\SS$- and $\SSSS$-partitions}\label{s:latin}

In this sections, we consider the equitable partitions of another graph,
$H(r,4)$, with the same quotient matrices \eqref{eq:S4}, \eqref{eq:S2}, 
and show how they are connected with $\SSSS$- and $\SS$-partitions.
The \emph{Hamming graph} $H(r,q)$ is the graph on the words of length $r$
over an alphabet $\Sigma$ of size $q$ (to be explicit, $\Sigma:=\{0,\ldots,q-1\}$),
two words being adjacent if and only if they differ in exactly one position.
For the purposes of our study,
we restrict consideration by the case $q=4$ only, including the definitions,
which can be easily expanded to an arbitrary $q$.
The graph $H(n,4)$ is regular of degree $3r$. 
A \emph{Latin $r$-cube} or \emph{Latin hypercube}, of order $4$ is an equitable
partition of $H(r,4)$ with the quotient matrix $\SSSS$, see \eqref{eq:S4}.
A set $C$ of vertices of $H(r,4)$  is a $4$-ary distance-$2$ \emph{MDS code}
(for brevity, we will omit ``$4$-ary distance-$2$'')
if $(C,\cC)$ is an equitable
partition of $H(r,4)$ with the quotient matrix $\SS$, see \eqref{eq:S2}.
The following known and straightforward correspondence is similar to 
Lemma~\ref{l:42}, but without any additional requirement on $C$.
\begin{lemma}\label{l:latinMDS} 
The Latin $r$-cubes $(C_0,C_1,C_2,C_3)$ of order $4$ 
are in one-to-one 
correspondence with the MDS codes $C$ in $H(r+1,4)$:
%\begin{IEEEeqnarray}{c}
%\label{eq:latinMDS}
$ C = \{ \vec v i \mid \vec v \in C_i, \ i=0,1,2,3\}.$
%\end{IEEEeqnarray}
\end{lemma}

As we see from Lemma~\ref{l:cov} below, the similarity between the quotient matrices corresponding to the 
Latin hypercubes on one side and the $\SSSS$-partitions on the other side is not an accident, and there is simple injective map from the first class of partitions to the second
(similarly, for MDS codes and $\SS$-partitions).
This connection is important because the class of Latin hypercubes of order $4$
was deeply studied in previous works and we can say that it is rather well understood.
In particular:
\begin{enumerate}
 \item There is a constructive characterisation of this class \cite{KroPot:4}:
 every Latin hypercube of order $4$ is semilinear or reducible. 
 The semilinear MDS codes correspond to semilinear $\SS$-partitions, 
 but in contrast with the class of all semilinear $\SS$-partitions,
 the analog of Theorem~\ref{th:semi} reduces the semilinear MDS codes
 (orthogonal to a fixed vector) to the class of Boolean functions, 
 which is trivially constructive.
 
 \item The asymptotic of the number $L(r)$ of objects is known \cite{PotKro:asymp}.
 In particular, the log\,log of this number is $\log_2\log_2 L(r)=r+o(r)$;
 we note that it corresponds to the simple lower bound obtained by a switching approach
 similar to \cite{Vas:nongroup_perfect.en} rather than the trivial upper bound $\log_2\log_2 L(r)<2r$.
 
 \item There are positive answers to the questions similar 
 to Problem~\ref{probl:2to4}
 and Problem~\ref{probl:3to4}. The first is almost trivial;
 the second has a complicate proof \cite{Pot:2012:partial},
 even using the characterisation \cite{KroPot:4}.
\end{enumerate}

All of those is applicable to the classes of $\SSSS$-partitions obtained from the 
Latin hypercubes.

Our main objects can be constructed from Latin hypercubes or MDS codes utilizing the following  
fact (for $S=\SS$, this is a partial case of \cite[Theorem~2]{Potapov:2010}). 
\begin{lemma}\label{l:cov}
 If $(C_i)_{i=0}^k$ is an equitable partition of $H(r,4)$ 
 with a quotient matrix $S$,
 then  $(C'_i)_{i=0}^k$ is an equitable partition of $Q_{3r}$
 with the quotient matrix $S$, where
\begin{IEEEeqnarray}{c} \label{eq:cov}
  C'_i := \!\!\!\!\!\! \bigcup_{(v_1,...,v_r)\in C_i} \!\!\!\!\!\! T_{v_1} \times \ldots \times T_{v_r},\quad  T_{0}:=\{000,111\},\\ \nonumber T_{1}:=\{010,101\},\ T_{2}:=\{100,011\},\ T_{3}:=\{110,001\}.
\end{IEEEeqnarray}
\end{lemma}
\begin{remark}
 In \eqref{eq:cov}, we see a correspondence between the vertices
 $(v_1,...,v_r)$ of $H(r,4)$ 
 and vertex subsets of $Q_{3r}$. This correspondence (to be exact, its inverse) 
 can be treated as a graph covering of $H(r,4)$ by $Q_{3r}$. In general,
 a  \emph{covering} of a graph $\mathrm{Tar}$ (the \emph{target})
 by the graph $\mathrm{Cov}$ (the \emph{cover})
 is an equitable partition of $\mathrm{Cov}$
 with the quotient matrix coinciding with the adjacency matrix of $\mathrm{Tar}$. 
 In the case of such covering, any equitable partition
 of $\mathrm{Tar}$ naturally induces an equitable partition of $\mathrm{Cov}$,
 called \emph{lifted} (this is a folklore fact).
 In particular, for any integer $q,m\ge 2$, there is an additive (over $\mathbb{Z}_q$)
 covering of $H(q^m,n)$ by $H(q,n(q^m-1)/(q-1))$.
\end{remark}

The previous lemma can also be treated as a kind of concatenated construction
for equitable partitions. In the case when $S$ is one of the quotient matrices considered in 
the current paper, the construction can be generalized, allowing to concatenate
equitable partitions in different dimensions and giving a reach classes
of $\SS$- and $\SSSS$-partitions including the ones that cannot be represented 
as in Lemma~\ref{l:cov} or Theorem~\ref{th:semi}.

\begin{theorem}\label{th:conc}
Assume that for each $j$ from $1$ to $r$,
$(T_i^{(j)})_{i=0,1,2,3}$ is an $\SSSS[n_i]$-partition, for some $n_i$.
 If $(C_i)_{i=0}^k$ is an equitable partition of $H(r,4)$ 
 with a quotient matrix $S^{(2)}_{3r}$, $S^{(3)}_{3r}$, or $S^{(4)}_{3r}$
 (where $k$ is $2$, $3$, or $4$, respectively),
 then  $(C'_i)_{i=0}^k$ is an equitable partition of $Q_{N}$,
 $N:=n_1+...+n_r$,
 with the quotient matrix $S^{(2)}_{N}$, $S^{(3)}_{N}$, or $S^{(4)}_{N}$, respectively.
\begin{IEEEeqnarray}{c} \label{eq:conc}
  C'_i := \!\!\!\!\!\! \bigcup_{(v_1,...,v_r)\in C_i} \!\!\!\!\!\! T^{(1)}_{v_1} \times \ldots \times T^{(r)}_{v_r}.  
\end{IEEEeqnarray}
\end{theorem}
\begin{IEEEproof}[A sketch of proof]
 The proof is straightforward, but utilizes the following fact,
 which is often used for equivalent definitions of MDS codes and Latin hypercubes:
 every cell $C'$ of an equitable partition of $H(r,4)$ with one of the considered quotient matrices
 intersects every clique of size $4$ in exactly $|C'|/4^{r-1}$ vertices. This explains why 
 the lemma cannot be generalized to an arbitrary quotient matrix.
\end{IEEEproof}

The $\SS$- and $\SSSS$-partitions that can be represented as in Theorem~\ref{th:conc} 
with $r\ge 0$ will be called \emph{reducible}.
In the remaining part of the paper, we discuss existence of  
$\SS$- and $\SSSS$-partitions that are neither semilinear nor reducible.

%========================================
%========================================
%========================================
\section{Irreducible full-rank partitions}\label{s:contr}
Theorems~\ref{th:semi} and~\ref{th:conc} give two powerful ways to construct 
$\SS$- and $\SSSS$-partitions (and hence, $(2n/3-1)$-resilient $(n,2)$-functions);
moreover, 
the concatenation construction of Theorem~\ref{th:conc} allows to combine
concatenated, or semilinear, or any other $\SSSS[n_i]$-partitions constructed in
smaller dimensions $n_i$.
Any partition that can be constructed using this two 
ways is semilinear or reducible (or both),
and it would be very important to find $\SS$- and $\SSSS$-partitions that are neither
semilinear nor reducible. 
Since at this moment no general method to construct such partitions is known,
we can only claim the following computational result.
\begin{proposition}\label{p:N16}
 There are an $\SS[12]$-partition and an $\SSSS[12]$-partition that are neither
semilinear nor reducible. 
\end{proposition}
The only $\SS[12]$-partition $(C,\cC)$ with this property is the last partition in
\cite[Appendix]{KroVor:CorIm}. There are at least $3$ disjoint translations
$C$, $C_1$, $C_2$ of $C$, providing a $\SSSS[12]$-partition 
$(C,C_1,C_2,\overline{C \cup C_1 \cup C_2})$ with the required property.

\begin{problem}
 Do there exist $\SS$- and $\SSSS$-partitions that are full rank and irreducible
 (i.e., neither semilinear nor reducible), for any $n\ge 12$, $n\equiv 0\bmod 3$?
\end{problem}

%========================================
%========================================
%========================================
\section{Conclusion}\label{s:conc}
\enlargethispage{-7cm} 
In this thesis, we discussed connections of 
the $(2n/3-1)$-resilient $(n,2)$-functions
(which has the best correlation-immunity order $(2n/3-1)$ 
among the balanced $(n,2)$-functions)
with special classes of equitable partitions, 
including the widely studied classes of $1$-perfect codes and Latin hypercubes.

The $(2n/3-1)$-resilient $(n,2)$-functions are equivalent to 
the equitable partition of a special type, called the $\SSSS$-partitions.
We described the non-full-rank $\SSSS$-partitions in terms
of multifold binary codes (Theorem~\ref{th:semi}), which gives a construction of at least 
$2^{2^{n/3}}$ different $\SSSS$-partitions.
Another presented construction (Theorem~\ref{th:conc}) 
uses $\SSSS[n_i]$-partitions with $n_i<n$ 
and Latin hypercubes of order $4$; we call resulting partitions \emph{reducible}.
No general way is known to construct full-rank irreducible 
$\SSSS$-partitions, but computational results say than in the dimension $n=12$
such partitions exist.
We conclude that the current knowledge is not sufficient to make 
a conjecture about the structure of an arbitrary $(2n/3-1)$-resilient $(n,2)$-function
(in contrast to the known characterization of the Latin hypercubes of order $4$),
but enough to construct double-exponential number 
of such functions with different properties.

\medskip
This work was funded by the Russian Science Foundation under grant 18-11-00136.

% 
% \bibliographystyle{IEEEtranS}
% \bibliography{../../k}

\begin{thebibliography}{10}
\providecommand{\url}[1]{#1}
\csname url@samestyle\endcsname
\providecommand{\newblock}{\relax}
\providecommand{\bibinfo}[2]{#2}
\providecommand{\BIBentrySTDinterwordspacing}{\spaceskip=0pt\relax}
\providecommand{\BIBentryALTinterwordstretchfactor}{4}
\providecommand{\BIBentryALTinterwordspacing}{\spaceskip=\fontdimen2\font plus
\BIBentryALTinterwordstretchfactor\fontdimen3\font minus
  \fontdimen4\font\relax}
\providecommand{\BIBforeignlanguage}[2]{{%
\expandafter\ifx\csname l@#1\endcsname\relax
\typeout{** WARNING: IEEEtranS.bst: No hyphenation pattern has been}%
\typeout{** loaded for the language `#1'. Using the pattern for}%
\typeout{** the default language instead.}%
\else
\language=\csname l@#1\endcsname
\fi
#2}}
\providecommand{\BIBdecl}{\relax}
\BIBdecl

\bibitem{Carlet:vect}
C.~Carlet, \emph{Vectorial Boolean Functions for Cryptography}, ser. Encycl.
  Math. Appl.\hskip 1em plus 0.5em minus 0.4em\relax Cambridge Univ. Press,
  2010, vol. 134, ch.~9, pp. 398--469, \DOI{10.1017/CBO9780511780448.012}.

\bibitem{FDF:CorrImmBound}
\BIBentryALTinterwordspacing
D.~G. Fon-Der-Flaass, ``A bound on correlation immunity,''
  \emph{\href{http://semr.math.nsc.ru}{Sib. Ehlektron. Mat. Izv.}}, vol.~4, pp.
  133--135, 2007, online: \url{http://mi.mathnet.ru/eng/semr149}. [Online].
  Available: \url{http://mi.mathnet.ru/eng/semr149}
\BIBentrySTDinterwordspacing

\bibitem{FDF:12cube.en}
\BIBentryALTinterwordspacing
------, ``Perfect colorings of the $12$-cube that attain the bound on
  correlation immunity,'' \emph{\href{http://semr.math.nsc.ru}{Sib. Ehlektron.
  Mat. Izv.}}, vol.~4, pp. 292--295, 2007, in Russian. English translation:
  \url{https://arxiv.org/abs/1403.8091}. [Online]. Available:
  \url{http://mi.mathnet.ru/eng/semr158}
\BIBentrySTDinterwordspacing

\bibitem{Friedman:92}
J.~Friedman, ``On the bit extraction problem,'' in \emph{Foundations of
  Computer Science, IEEE Annual Symposium on}.\hskip 1em plus 0.5em minus
  0.4em\relax Los Alamitos, CA, USA: IEEE Computer Society, 1992, pp. 314--319,
  \DOI{10.1109/SFCS.1992.267760}.

\bibitem{KroPot:nonsplittable}
D.~S. Krotov and V.~N. Potapov, ``On multifold {MDS} and perfect codes that are
  not splittable into onefold codes,''
  \emph{\href{http://link.springer.com/journal/11122}{Probl. Inf. Transm.}},
  vol.~40, no.~1, pp. 5--12, 2004, \DOI{10.1023/B:PRIT.0000024875.79605.fc}
  translated from
  \href{http://www.mathnet.ru/php/journal.phtml?jrnid=ppi\&option_lang=eng}{Probl.
  Peredachi Inf.} 40(1) (2004), 6-14.

\bibitem{KroPot:4}
------, ``$n$-{A}ry quasigroups of order $4$,''
  \emph{\href{http://epubs.siam.org/journal/sjdmec}{SIAM J. Discrete Math.}},
  vol.~23, no.~2, pp. 561--570, 2009, \DOI{10.1137/070697331}.

\bibitem{KroVor:CorIm}
\BIBentryALTinterwordspacing
D.~S. Krotov and V.~K. Vorob'ev, ``On unbalanced {B}oolean functions attaining
  the bound $2n/3-1$ on the correlation immunity,'' arXiv.org, E-print
  1812.02166v2, 2018. [Online]. Available:
  \url{https://arxiv.org/abs/1812.02166}
\BIBentrySTDinterwordspacing

\bibitem{Potapov:2010}
\BIBentryALTinterwordspacing
V.~N. Potapov, ``On perfect colorings of {B}oolean $n$-cube and correlation
  immune functions with small density,''
  \emph{\href{http://semr.math.nsc.ru}{Sib. Ehlektron. Mat. Izv.}}, vol.~7, pp.
  372--382, 2010, in Russian, English abstract. [Online]. Available:
  \url{http://mi.mathnet.ru/eng/semr248}
\BIBentrySTDinterwordspacing

\bibitem{Pot:2012:partial}
------, ``On extensions of partial $n$-quasigroups of order $4$,''
  \emph{\href{http://link.springer.com/journal/12002}{Sib. Adv. Math.}},
  vol.~22, no.~2, pp. 135--151, 2012, \DOI{10.3103/S1055134412020058}
  translated from
  \href{http://www.mathnet.ru/php/journal.phtml?jrnid=mt\&option_lang=eng}{Mat.
  Tr.} 14(2):147-172, 2011.

\bibitem{Pot:2012:color}
------, ``On perfect $2$-colorings of the $q$-ary $n$-cube,''
  \emph{\href{http://www.sciencedirect.com/science/journal/0012365X}{Discrete
  Math.}}, vol. 312, no.~6, pp. 1269--1272, 2012,
  \DOI{10.1016/j.disc.2011.12.004}.

\bibitem{PotKro:asymp}
V.~N. Potapov and D.~S. Krotov, ``Asymptotics for the number of $n$-quasigroups
  of order $4$,'' \emph{\href{http://link.springer.com/journal/11202}{Sib.
  Math. J.}}, vol.~47, no.~4, pp. 720--731, 2006,
  \DOI{10.1007/s11202-006-0083-9} translated from
  \href{http://www.mathnet.ru/php/journal.phtml?jrnid=smj\&option_lang=eng}{Sib.
  Mat. Zh.} 47(4) (2006), 873-887.

\bibitem{Tarannikov2000}
\BIBentryALTinterwordspacing
Y.~Tarannikov, ``On resilient {B}oolean functions with maximal possible
  nonlinearity,'' Cryptology ePrint Archive 2000/005, 2000,
  \url{https://eprint.iacr.org/2000/005}. [Online]. Available:
  \url{https://eprint.iacr.org/2000/005}
\BIBentrySTDinterwordspacing

\bibitem{Vas:nongroup_perfect.en}
Y.~L. Vasil'ev, ``On nongroup close-packed codes,'' in \emph{Probleme der
  Kybernetik}.\hskip 1em plus 0.5em minus 0.4em\relax Akademie-Verlag, 1965,
  vol.~8, pp. 92--95, translated from Problemy Kibernetiki 8: 337-339, 1962.

\end{thebibliography}
% \end{document}

% Generated by IEEEtranS.bst, version: 1.14 (2015/08/26)
\providecommand\href[2]{#2} \providecommand\url[1]{\href{#1}{#1}}
  \def\DOI#1{{\small {DOI}:
  \href{http://dx.doi.org/#1}{#1}}}\def\DOIURL#1#2{{\small{DOI}:
  \href{http://dx.doi.org/#2}{#1}}}

\end{document}